\documentclass{amsart}
\newtheorem{thm}{Theorem}
\newtheorem{note}{Note}
\newtheorem{prop}{Proposition}
\begin{document}

\title{About coordinates on the phase-spaces of\\
Schlesinger system ($n+1$ matrices, $sl(2,\mathbb{C})$-case)\\ and Garnier--Painlev\'e~6 system}%
\address{PDMI, Fontanka 27, St.Petersburg, 191023, Russia}%
\email{mbabich@pdmi.ras.ru}%
\author{Mikhail V. Babich}
\keywords{Schlesinger system, Painlev\'e equations, isomonodromic deformations, phase space, simplectic coordinates}%
\begin{abstract}
The geometric model of the pathway linking Schlesinger and
Garnier--Painlev\'e~6 systems based on an original orthonormalization of a set
of elements in ${sl(2,\mathbb C)}$ is constructed. The explicit
polynomial map of the Cartesian products of $n-2$ quadrics (the
Zariski-topology chart of the phase space of the
Garnier--Painlev\'e~6 system) into the phase space of the Schlesinger
system and  the rational inverse to this map are presented.
\end{abstract}
\maketitle

The subject of our considerations will be the phase space of the
Schlesinger system (SchS) of equations (about SchS and around it see, for example,
\cite{GtoP}):
\begin{equation}\label{eq1}
dA^{(k)}=\sum_{i\ne k}
[A^{(i)},A^{(k)}]d\log (\lambda_i-\lambda_k), \ \ \
\sum_k A^{(k)}=0
\end{equation}
for $n+1$ ${sl(2,\mathbb C)}$-valued matrices $A^{(k)}$ depending on complex
parameters $\lambda_0, \lambda_1, \dots , \lambda_n$.

Let us introduce the notations. We denote $\langle A,B\rangle =tr \
AB$ --- the Killing product. It set the isomorphism between ${sl(2,\mathbb C)}$
and $sl^*(2,\mathbb C)$;
by the matrix elements of $|A\rangle \in sl^*(2,\mathbb C)$ we mean
the matrix elements of the corresponding $A\in {sl(2,\mathbb C)}$.

Denote
$a_{ij}:=\langle A^{(i)},A^{(j)}\rangle, \ \  f_{ijk}:=\langle [A^{(i)},A^{(j)}],A^{(k)}\rangle$.
The values $a_{ii}=-2\det A^{(i)}$, $i=0,1,\dots ,n$ keep constant
values for the solutions of (\ref{eq1}), they are the parameters of
SchS.

Let $SL|\stackrel{\scriptscriptstyle a_{kk}}{A}\rangle$ be the orbit
of the co-adjoint action of $SL(2, \mathbb{C})$
on nonzero element $|\stackrel{\scriptscriptstyle a_{kk}}{A}\rangle$ such
that $\langle \stackrel{\scriptscriptstyle
a_{kk}}{A},\stackrel{\scriptscriptstyle a_{kk}}{A}\rangle=a_{kk}$. Note that
$A^{(k)}\in SL|\stackrel{\scriptscriptstyle
a_{kk}}{A}\rangle $ if $A^{(k)}\ne 0$ only,
so for zero values of $A^{(k)}$ we need to extend the orbit.
We blow up the vertex of cone $\langle \stackrel{\scriptscriptstyle
0}{A},\stackrel{\scriptscriptstyle 0}{A}\rangle=0$.

We define $SL|{A}\rangle ' $, $\langle A,A\rangle=R^2$ as the
submanifold of ${sl(2,\mathbb C)} \times \mathbf{P}SL|A\rangle$.
Its  points have ``the affine component''
$A$ and ``the projective component'' $\tilde A$:
$$\left(A,\tilde A \right)=\left(
\left(%
\begin{array}{cc}
  X_3 & X_1 \\
  X_2 & -X_3 \\
\end{array}%
\right) ,
\left(%
\begin{array}{cc}
  x_3 & x_1 \\
  x_2 & -x_3 \\
\end{array}%
\right)\right)  \in SL|{A}\rangle '.$$
The submanifold is defined by the equations
$$
X_1X_2+X_3^2=R^2/2, \  \
x_1x_2+x_3^2=x_0^2R^2/2, \ \
X_i=x_0x_i, \ \ i= 1,2,3.
$$

We can see that $SL|{A}\rangle ' $ is an algebraic simplectic manifold with the form
induced by the restriction of the Lie-Poisson bracket on the orbit:
$$\omega_{\scriptscriptstyle LP}^R = dX_3\wedge d\log X_1=-dX_3\wedge d\log X_2
=d\frac{X_2}{\sqrt{R^2/2}-X_3}\wedge dX_1=
d\frac{X_1}{\sqrt{R^2/2}+X_3}\wedge dX_2\ ;$$
if $R^2=0$,  $\omega_{\scriptscriptstyle LP}^0 =d(x_3:x_1)\wedge dX_1=d(-x_3:x_2)\wedge dX_2$.

 The simplectic manifold $SL|{A}\rangle ' $
can be covered by two standard simplectic charts
$(\mathbb{A}^2, dp\wedge dq)$, $(\mathbb{A}^2, dp'\wedge dq')$:
with the transition functions
$p'=1/p, \ \ \ q'=-q(pq+2\sqrt{R^2/2})$, for example
$(p,q)=((X_3-\sqrt{R^2/2})/{X_1},X_1)$,
$(p',q')=({X_1}/(X_3-\sqrt{R^2/2}),X_2)$.

Traditionally the manifold $SL|{A}\rangle ' $ is defined as the abstract manifold,
covered by these two charts, see \cite{Ok-s,ShT}.

\begin{note}
The manifold $SL|A\rangle '$ is simplectomorphic
to the smooth quadric $X_1X_2+X_3^2=R^2/2$,  $R^2\ne 0$
or to the cone with blowing up  vertex;
these simplectomorphisms are rational.
\end{note}

We call $SL|\stackrel{a_{kk}}{A}\rangle '$ ``the quadric'' in all cases, there are the quadrics that
we talk about in the Abstract,
\emph{the phase space of the SchS is the submanifold $\sum A^{(i)}=0$ of the Cartesian product of such quadrics.}

\underline{\textbf{Restriction}}

\emph{We exclude the special case, the sets all
the matrices $A^{(i)}$, $i=0, \dots, n$ can be
carried into the upper-triangle form simultaneously.}

This case is much more simple then
the general one, but needs another method;
let us denote by
$\Delta_{a_{\overline{ii}}}$ such ``all-triangle'' sets $A^{(\overline{i})}$ and put

$$\stackrel{\approx}{M}_{a_{\overline{ii}}}:=
\left( SL|\stackrel{\scriptstyle a_{00}}{A}\rangle '\times
SL|\stackrel{\scriptstyle a_{11}}{A}\rangle '\times \dots \times
SL|\stackrel{\scriptstyle a_{nn}}{A}\rangle '\right)
\setminus \Delta_{a_{\overline{ii}}}
$$

The bar over index means it is the
set of such values with all values of indexes:
$\stackrel{\approx}{M}_{a_{\overline{ii}}}=\stackrel{\approx}{M}_{a_{00},a_{11}, \dots, a_{nn}}$;
we will use the similar notation  --
$A^{(\overline{i})}$ is equivalent to $A^{(0)},A^{(1)},
\dots , A^{(n)}$ etc.

Manifold $\stackrel{\approx}{M}_{a_{\overline{ii}}}$ is $2(n+1)$-dimensional simplectic space with the form
$
    \stackrel{\approx}{\omega}:=\sum_{k=0}^n
 \omega^{(k)},
 $
 ${\omega}^{(k)}$ is the form $\omega^{\sqrt{a_{kk}}}_{\scriptscriptstyle LP}$ on the ``$k$''-th
 Cartesian factor $SL|\stackrel{\scriptstyle a_{kk}}{A}\rangle '$.

The group ${SL(2,\mathbb C)}$ acts on $\mathbf{P}SL|A\rangle$
(the projectivization commutes with the co-adjoint action).
We define (the component-wise) action of ${SL(2,\mathbb C)}$ on
$\stackrel{\approx}{M}_{a_{\overline{ii}}}$ and
denote this action by a sub-index:
$$g\left({A'}^{(\overline{i})}\right):={A'}^{(\overline{i})}_g:=
(A^{(0)}_g,{\tilde A}^{(0)}_g),\dots ,
(A^{(n)}_g,{\tilde A}^{(n)}_g), \ \ A^{(i)}_g:=gA^{(i)}g^{-1}, \ {\tilde A}^{(i)}_g:=
g{\tilde A}^{(i)}g^{-1}.
$$

Denote the submanifold of $\stackrel{\approx}{M}_{a_{\overline{ii}}}$ that
consists of such ${A'}^{(\overline{i})}$ that $\sum A^{(k)}=0$ by
$\stackrel{\sim}{M}_{a_{\overline{ii}}}$ .
It is evident that such action of ${SL(2,\mathbb C)}$ preserve the property $\sum A^{(k)}=0$,
so group ${SL(2,\mathbb C)}$ acts on $\stackrel{\sim}{M}_{a_{\overline{ii}}}$ too.

\emph{The phase space $M_{a_{\overline{ii}}}$ of the Garnier system (if $n=3$ it is
the Painlev\'e~6-system) is the quotient of $\stackrel{\sim}{M}_{a_{\overline{ii}}}$ with respect
to this action:
$M_{a_{\overline{ii}}}:=\stackrel{\sim}{M}_{a_{\overline{ii}}}/{SL(2,\mathbb C)}$}, it is the main
object of our building.

We need functions (coordinates) on the quotient of with respect to ${SL(2,\mathbb C)}$.
The method is based on the proposition.

Let $\sigma^
{(\overline{i})}=\sigma^{(1)},\sigma^{(2)},\sigma^{(3)}$ be a basis
of ${sl(2,\mathbb C)}$, let
$\sigma^{(\overline{i})}_g=g\sigma^{(1)}g^{-1},g\sigma^{(2)}g^{-1},
g\sigma^{(3)}g^{-1}$ be this basis ``turned'' by an element $g \in
{SL(2,\mathbb C)}$, and $SL\sigma^{(\overline{i})}:=\bigcup_{g\in
{SL(2,\mathbb C)}}\sigma^{(\overline{i})}_g$ be the orbit of the
action of ${SL(2,\mathbb C)}$ on $\sigma^{(\overline{i})}$.

Let the map $\stackrel{\sim}{M}_{a_{\overline{ii}}}\to
SL\sigma^{(\overline{i})}$: ${A'}^{(\overline{i})}\to
\sigma^{(1)}({A'}^{(\overline{i})}),\sigma^{(2)}({A'}^{(\overline{i})}),
\sigma^{(3)}({A'}^{(\overline{i})})
$
commutes with the action of ${SL(2,\mathbb C)}$:
$$\sigma^{(1)}({A'}_g^{(\overline{i})}),\sigma^{(2)}({A'}_g^{(\overline{i})}),
\sigma^{(3)}({A'}_g^{(\overline{i})})=
\sigma_g^{(1)}({A'}^{(\overline{i})}),\sigma_g^{(2)}({A'}^{(\overline{i})}),
\sigma_g^{(3)}({A'}^{(\overline{i})}).
$$

We call this property \emph{the ${SL(2,\mathbb C)}$-invariance}.
\begin{prop}\label{proposition1}
The coordinates of all $A^{(k)}$ in the basis
$\sigma^{(\overline{i})}({A'}^{(\overline{i})})$ constructed by an ${SL(2,\mathbb C)}$-invariant map
do not depend on the
action of ${SL(2,\mathbb C)}$ on $\stackrel{\sim}{M}_{a_{\overline{ii}}}$ , they are functions
on $\stackrel{\sim}{M}_{a_{\overline{ii}}}/{SL(2,\mathbb C)} $.
$\Box$
\end{prop}

There are the functions we will use as the coordinates on ${M}_{a_{\overline{ii}}}$ .

\begin{note}
A related map is the routine orthonormalization of a basis in  ${\mathbb R}^n$,
but we will use some special method based on the Lie-structure of ${sl(2,\mathbb C)}$.
\end{note}

Consider such bases $\sigma'_{\pm,3}$ that
\begin{equation}\label{eq2}
\langle \sigma'_{-},\sigma'_{-}\rangle
=\langle \sigma'_{+},\sigma'_{+}\rangle
=\langle \sigma'_{-},\sigma'_{3}\rangle
=\langle \sigma'_{+},\sigma'_{3}\rangle =0, \ \
\langle \sigma'_{+},\sigma'_{-}\rangle=1, \ \ [\sigma'_{+},\sigma'_{-}]=\sigma_3.
\end{equation}
Denote the set of these bases $SL\sigma_{\pm,3}$, we call them
``\emph{the standard}'' bases. The foundation of the offered method
is the following proposition.

\begin{prop}
For any $A^{(0)}\ne 0$ and any sign of the square root, the value
$\sqrt{2a_{00}}:=\sqrt{-4 \det A^{(0)}}$
is the eigenvalue of the linear operator
$ad_{A^{(0)}}=[A^{(0)}, \cdot ]: \ \ {sl(2,\mathbb C)}\to{sl(2,\mathbb C)}$.
The corresponding eigenspace is one-dimensional, isotropic and orthogonal to $A^{(0)}$.
\end{prop}
\textbf{Proof.}
All statements are ${SL(2,\mathbb C)}$-invariant, consequently we can carry $A^{(0)}$ into the Jordan form.
For such matrices the Proposition is evident.
$\Box$

Denote this eigenspace by $(\sigma^{(\sqrt{2a_{00}})})$, and by $\sigma^{(\sqrt{2a_{00}})}$
any eigenvector from it.
We define $(\sigma^{(\sqrt{2a_{ii}})})$ for ${A'}^{(i)}$
from the blowing up vertices now. If $a_{ii}=0$ and
${A}^{(i)}\ne 0$, $(\sigma^{(\sqrt{2a_{ii}})})$ is the direction of ${\tilde A}^{(i)}$,
we set \emph{
$(\sigma^{(\sqrt{2a_{ii}})})$ the direction of ${\tilde A}^{(i)}$ for ${A}^{(i)}= 0$ too.}

Consider three vectors $A^{(0)}$, $A^{(n-1)}$, $A^{(n)}$ and fix any values of
$\sqrt{2a_{00}}$, $\sqrt{2a_{n-1n-1}}$.

We call the standard basis \emph{accompanying to $A^{(0)}$, $A^{(n-1)}$, $A^{(n)}$ along
$(\sigma^{(-\sqrt{2a_{00}})})$ and $(\sigma^{(\sqrt{2a_{n-1n-1}})})$} if
\begin{equation}\label{eq3}
    \sigma_{-}\in (\sigma^{(-\sqrt{2a_{00}})}), \ \
    \sigma_{+}\in (\sigma^{(\sqrt{2a_{n-1n-1}})}), \ \
    \langle \sigma_{-}, A^{(n)}\rangle =1.
\end{equation}
\begin{prop}
The standard basis accompanying to ${A'}^{(0)}$, ${A'}^{(n-1)}$, ${A'}^{(n)}$ along
$(\sigma^{(-\sqrt{2a_{00}})})$ and $(\sigma^{(\sqrt{2a_{n-1n-1}})})$ exists and unique
if and only if
\begin{equation}\label{eq4}
      \text{\textbf{a)}} \ \ \langle \sigma^{(-\sqrt{2a_{00}})}, \sigma^{(\sqrt{2a_{n-1n-1}})}\rangle  \ne  0
\text{ , and \ \ \ }      \text{\textbf{b)}}  \ \ \langle \sigma^{(-\sqrt{2a_{00}})}, A^{(n)}\rangle  \ne  0.
\end{equation}
\end{prop}
\textbf{Proof.}
If ``\textbf{a)}'' and ``\textbf{b)}'' are satisfied, the explicit formulae
\begin{equation}\label{eq5}
    \begin{array}{ccl}
      \sigma_{-} & = & \sigma^{(-\sqrt{2a_{00}})}
      \frac{1}{\langle \sigma^{(-\sqrt{2a_{00}})}, A^{(n)}\rangle} \\
      \sigma_{+} & = & \sigma^{(\sqrt{2a_{n-1n-1}})}
      \frac{\langle \sigma^{(-\sqrt{2a_{00}})}, A^{(n)}\rangle}{
      \langle \sigma^{(-\sqrt{2a_{00}})},\sigma^{(\sqrt{2a_{n-1n-1}})} \rangle} \\
      \sigma_{3} & = &  [\sigma^{(\sqrt{2a_{n-1n-1}})}, \sigma^{(-\sqrt{2a_{00}})} ]
      \frac{1}{\langle \sigma^{(-\sqrt{2a_{00}})},\sigma^{(\sqrt{2a_{n-1n-1}})} \rangle}
          \end{array}
\end{equation}
give the desired basis. It is unique.

If ``\textbf{a)}'' does not satisfied, vectors $\sigma_{-}$ and $\sigma_{+}$ are linear dependent,
$\sigma_{\pm,3}$ is not a basis. If ``\textbf{b)}'' does not satisfied, (\ref{eq3}) can not take place.
$\Box$

The subset of $\stackrel{\sim}{M}_{a_{\overline{ii}}}$ where (\ref{eq4}) fulfilled we
denote by $\stackrel{\sim}{\mathcal{U}}(-\sqrt{2a_{00}}//\sqrt{2a_{n-1n-1}}/n)$, or $\stackrel{\sim}{\mathcal{U}}(\ \dots)$
for short:
$$\stackrel{\sim}{\mathcal{U}}(\ \dots):=\{{A'}^{(\overline{i})}\in \stackrel{\sim}{M}_{a_{\overline{ii}}}\ : \ \
\langle \sigma^{(-\sqrt{2a_{00}})}, \sigma^{(\sqrt{2a_{n-1n-1}})}\rangle  \ne  0 , \ \
 \langle \sigma^{(-\sqrt{2a_{00}})}, A^{(n)}\rangle  \ne 0
 \}
$$

From the Restriction follows
\begin{thm}
For any preassigned values of square roots $\sqrt{2a_{ii}}, \ \ i\in \{0,1,\dots,n\}$
$$\stackrel{\sim}{M}_{a_{\overline{ii}}}=
\bigcup_{\left. i,j\atop 0\ne i \ne j\ne 0\right. }\stackrel{\sim}{\mathcal{U}}(-\sqrt{2a_{00}}//\sqrt{2a_{ii}}/j)
$$
\end{thm}
$\Box$

The Restriction implies that the values $a_{ik}$, $f_{ijk}$ characterize the set $A^{(\overline{i})}$
up to the action of $SL(2,\mathbb C)$, consequently the quotient
$M_{a_{\overline{ii}}}:=\stackrel{\sim}{ M}_{a_{\overline{ii}}}/{SL(2,\mathbb C)}$
\emph{is the manifold }. If $a_{kk}\ne 0$ it can be embedded into the space of all $a_{ij}$, $f_{ijk}$;
for the sets with ${A'}^{(k_0)}$ on the blowing up divisor
some coordinates should be taken from the conic
$\langle {\tilde A}^{(k_0)}, {\tilde A}^{(k_0)}\rangle=0$ on the plane
$\mathbb{C}\mathbf{P}^2\ni \langle {\tilde A}^{(k_0)}, {A}^{(i_1)}\rangle:\langle {\tilde A}^{(k_0)},
{A}^{(i_2)}\rangle:\langle {\tilde A}^{(k_0)}, [{A}^{(i_1)},{A}^{(i_2)}]\rangle$.

We treat $M_{a_{\overline{ii}}}$ as the abstract manifold. The quotation
of $\stackrel{\sim}{ M}_{a_{\overline{ii}}}$ with respect to ${SL(2,\mathbb C)}$ we denote by $\tilde \pi$:
$$\tilde \pi: \ \ \stackrel{\sim}{ M}_{a_{\overline{ii}}} \longrightarrow
M_{a_{\overline{ii}}}:=\stackrel{\sim}{ M}_{a_{\overline{ii}}}/{SL(2,\mathbb C)}  .
$$

Denote the embedding of $M_{a_{\overline{ii}}}$ into the
space of values $a_{ij}$, $f_{ijk}$
 and ratios
$\langle A^{(i)}, {\tilde A}^{(k_0)}\rangle:\langle [A^{(i)}, A^{(j)}],{\tilde A}^{(k_0)}\rangle$,
by $M^f_{a_{\overline{ii}}}$.

Consider an another version of the quotation now. Formulae (\ref{eq5}) give the map
$A^{(\overline{k})} \longrightarrow \sigma_{\pm,3}(A^{(\overline{k})})$ and basis
that we talk about in the Proposition~\ref{proposition1}.
The map is defined on $\stackrel{\sim}{\mathcal{U} }(-\sqrt{2a_{00}}//\sqrt{2a_{n-1n-1}}/n)
\subset \stackrel{\sim}{M}_{a_{\overline{ii}}}$ and
induce  the embedding of
$\stackrel{\sim}{\mathcal{U} }(-\sqrt{2a_{00}}//\sqrt{2a_{n-1n-1}}/n)/{SL(2,\mathbb C)}
\subset {M}_{a_{\overline{ii}}}$
into $\stackrel{\approx}{M}_{a_{\overline{ii}}}$ \ :

\begin{equation}\label{eq6}
\begin{array}{ll}
A^{(0)}=\left(%
\begin{array}{cc}
  \sqrt{\frac{a_{00}}{2}} & 0 \\
  q'_0 & -\sqrt{\frac{a_{00}}{2}} \\
\end{array}%
\right) &
A^{(i)}=\left(%
\begin{array}{cc}
  \beta_i & q_i \\
  q'_i & -\beta_i \\
\end{array}%
\right)\\
A^{(n-1)}=\left(%
\begin{array}{cc}
  \sqrt{\frac{a_{n-1n-1}}{2}} & q_{n-1} \\
  0 & -\sqrt{\frac{a_{n-1n-1}}{2}} \\
\end{array}%
\right)  &
A^{(n)}=\left(%
\begin{array}{cc}
  \beta_n & 1 \\
  q'_n & -\beta_n \\
\end{array}%
\right)
\end{array}
\end{equation}
 -- it is the form the matrices from the set $A^{(\overline{i})}$ have in the
 accompanying basis $\sigma_{\pm,3}$.
This map is the embedding because in the fixed basis every vector (matrix) is uniquely
defined by its coordinates (matrix elements).

Let us reject $A^{(0)},A^{(n-1)},A^{(n)}$ from the set $A^{(\overline{i})}$, it is the projection
to the Cartesian product of all $SL|\stackrel{a_{ii}}{A}\rangle '$ except those with $i=0,n-1,n$
$$\stackrel{\approx}{M}_{a_{\overline{ii}}\setminus \{0n-1n\}}:= \ \ SL|\stackrel{a_{11}}{A}\rangle '  \times
 \dots \times SL|\stackrel{a_{n-2n-2}}{A}\rangle ' \ ,$$
it is the mentioned in the Abstract product of $n-2$ quadrics.

This projection is a bijection on the image,
the rejected terms can be restored. Denote
$$
\beta_\Sigma:=\sum\limits^{n-2}_{i=1}\beta_i, \ \
q_\Sigma:=\sum\limits^{n-2}_{i=1} q_i, \ \
q'_\Sigma:=\sum\limits^{n-2}_{i=1} q'_i \ .
$$
Because of $\sum A^{(i)}=0$ and $\det A^{(n)}=-a_{nn}/2$
\begin{equation}\label{eq7}
\begin{array}{ll}
&-q_{n-1}=q_{\Sigma}+1,\ q'_{n-1}:=0,\ \beta_{n-1}:=\sqrt{a_{n-1n-1}/2},\\
&-\beta_n=\sqrt{a_{00}/2}+\sqrt{a_{n-1n-1}/2}+\beta_\Sigma, \
q'_n=-(\beta_\Sigma+\sqrt{a_{00}/2}+\sqrt{a_{n-1n-1}/2})^2+a_{nn}/2,\\
&-q'_0=q'_\Sigma-(\beta_\Sigma+\sqrt{a_{00}/2}+\sqrt{a_{n-1n-1}/2})^2+a_{nn}/2,
\ \ q_0:=0, \ \ \beta_0:=\sqrt{a_{00}/2} \ .
\end{array}
\end{equation}

Denote the composition of the embedding (\ref{eq6}) and the rejection
of $A^{(0)},A^{(n-1)},A^{(n)}$ by $\pi_{\{0n-1n\}}$:
$$\pi_{\{0n-1n\}}: \ \ \stackrel{\sim}{\mathcal{U} }(-\sqrt{2a_{00}}//\sqrt{2a_{n-1n-1}}/n)/{SL(2,\mathbb C)}
\longrightarrow \stackrel{\approx}{M}_{a_{\overline{ii}}\setminus \{0n-1n\}}
$$
\begin{prop}
Map $\pi_{\{0n-1n\}}$ is the bijection. \ \ \ $\Box$
\end{prop}
The space $\stackrel{\approx}{M}_{a_{\overline{ii}}\setminus \{0n-1n\}}$ is the
simplectic manifold, as any product of simplectic spaces.
Denote its form by $\stackrel{\approx}{\omega}_{\{0n-1n\}}$.

\begin{thm}
The form $\stackrel{\approx}{\omega}|_{\stackrel{\sim}{M}\subset\stackrel{\approx}{M}}$,
that is the restriction of the simplectic form $\stackrel{\approx}{\omega}$
 on $\stackrel{\sim}{\mathcal{U} }(-\sqrt{2a_{00}}//\sqrt{2a_{n-1n-1}}/n) \subset
\stackrel{\sim}{M}_{a_{\overline{ii}}}\subset \stackrel{\approx}{M}_{a_{\overline{ii}}}$~,
coincides with the lifting of $\stackrel{\approx}{\omega}_{\{0n-1n\}}$ on
$\stackrel{\sim}{\mathcal{U} }(\ \dots) $:
${\tilde \pi}^{*}\circ {\pi}^{*}_{\{0n-1n\}}
\stackrel{\approx}{\omega}_{\{0n-1n\}}=
\stackrel{\approx}{\omega}|_{\stackrel{\sim}{M}\subset\stackrel{\approx}{M}}
$.
\end{thm}

\textbf{Proof.}
The restriction of $\stackrel{\approx}{\omega}$
on $\stackrel{\sim}{M}_{a_{\overline{ii}}}$ does not depend on the choice of basis of ${sl(2,\mathbb C)}$, so we
can calculate the sum
$\sum_{i=0}^n\omega^{(i)}$ in the accompanying along
$(\sigma^{(-\sqrt{2a_{00}})})$ and $(\sigma^{\sqrt{(2a_{n-1n-1})}})$ basis
that exists for ${A'}^{(\overline{i})}\in
\stackrel{\sim}{\mathcal{U} }(-\sqrt{2a_{00}}//\sqrt{2a_{n-1n-1}}/n) $.

The point is, in this basis
$\omega^{(0)}=\omega^{(n-1)}=\omega^{(n)}=0$.
It is true because ${A'}^{(0)}$, ${A'}^{(n-1)}$, ${A'}^{(n)}$   belong to
\emph{the one-dimensional} submanifolds of $SL|\stackrel{a_{ii}}{A}\rangle '$  -- to the
intersections of quadrics $SL|\stackrel{a_{ii}}{A}\rangle $, ($i=0,n-1,n$)
and planes $A^{(0)}_{12}=0$, $A^{(n-1)}_{21}=0$,
$A^{(n)}_{12}=1$, consequently
$\sum_{i=0}^n\omega^{(i)}=\sum_{i=1}^{n-2}\omega^{(i)}$ . \ \ $\Box$

Denote
$
\mathcal{U}(\ \dots ):=\tilde \pi
(\stackrel{\sim}{\mathcal{U}}(\ \dots ))\subset M_{a_{\overline{ii}}} \ .
$
The definition of
 $\stackrel{\sim}{\mathcal{U}}(\ \dots )$
is ${SL(2,\mathbb C)}$-invariant, consequently
$ \stackrel{\sim}{\mathcal{U}}(\ \dots )=
 \tilde{\pi}^{-1}({\mathcal{U}}(\ \dots ))$,
$$M_{a_{\overline{ii}}} =
\bigcup_{\left. i,j\atop 0\ne i \ne j\ne 0\right. }{\mathcal{U}}(-\sqrt{2a_{00}}//\sqrt{2a_{ii}}/j) \ ,
$$
and manifold \emph{$M_{a_{\overline{ii}}}$ is simplectic manifold}.
There is the global simplectic form $\omega$ can be glued from the forms on
${\mathcal{U}}(\ \dots )\simeq \stackrel{\approx}{M}_{a_{\overline{ii}}\setminus \{0ij\}}$
\begin{note}
The set of local simplectic coordinates on $M_{a_{\overline{ii}}}$
is a set $(p_i,q_i)_{i=1}^{n-2}$, where $(p_i,q_i)$ is any pair of local simplectic coordinates on the
quadric $SL|\stackrel{a_{ii}}{A}\rangle ' $.
\end{note}
\begin{thm}
Coordinate functions $a_{ik},f_{ijk}$ and the functions  $\beta_i,q_i,q'_i$,
the matrix elements of $A^{(i)}$ in the accompanying basis,
are birationally connected.
\end{thm}

\textbf{Proof.}
In one direction it is trivial, using the representation (\ref{eq6}) and
formulae (\ref{eq7}) we calculate $tr \ A^{(i)}A^{(j)}=a_{ij} $
and $tr \ [A^{(i)},A^{(j)}]A^{(k)}=f_{ijk} $. They will be some polynomials of matrix elements
$\beta_i, q_i, q'_i, \ \ i=1, \dots, n-2$ .

Consider the opposite  direction. The foundation of the construction is the following proposition
that can be verified by direct calculation:
\begin{prop}
For any $A,B\in {sl(2,\mathbb C)}$ vector $\sigma=\sigma^{(\sqrt{2\langle A,A\rangle})}(B\setminus A)$:
\begin{equation}\label{eq8}
\sigma^{(\sqrt{2\langle A,A\rangle})}(B\setminus A):=
\langle A,A\rangle B- \langle A,B\rangle A+\sqrt{\langle A,A\rangle /2}[A,B]
\end{equation}
satisfy the equality $[A,\sigma]=\sqrt{2\langle A,A\rangle} \sigma$. $\Box$
\end{prop}
The Restriction guarantee for any set $A^{(\overline{i})}$  there are such
$A^{\widehat{0}}$ and $A^{\widehat{n-1}}$ that
$\sigma^{(-\sqrt{2a_{00}})}(A^{(\widehat{0})}\setminus A^{(0)})\ne 0$,
$\sigma^{(\sqrt{2a_{n-1n-1}})}(A^{(\widehat{n-1})}\setminus A^{(n-1)})\ne 0$.
We set
\begin{equation}\label{eq9}
\sigma^{(-\sqrt{2a_{00}})}=\sigma^{(-\sqrt{2a_{00}})}(A^{(\widehat{0})}\setminus A^{(0)}), \ \ \ \
\sigma^{(\sqrt{2a_{n-1n-1}})}:=\sigma^{(\sqrt{2a_{n-1n-1}})}(A^{(\widehat{n-1})}\setminus A^{(n-1)}).
\end{equation}

For $A^{(\overline{i})}\in \stackrel{\sim}{\mathcal{U}}(-\sqrt{2a_{00}}//\sqrt{2a_{n-1n-1}}/n)$ vectors
$ \sigma^{(-\sqrt{2a_{00}})}$ and $\sigma^{(\sqrt{2a_{n-1n-1}})}$ are linear independent and
(\ref{eq5}) for (\ref{eq9}) give the standard basis $\sigma_{\pm,3}$ accompanying $A^{(\overline{i})}$.
The matrix elements $\beta_{i},q_i,q'_i$ are
values $\langle \sigma_{3}/2,A^{(i)}\rangle$,
$\langle \sigma_{-},A^{(i)}\rangle$, $\langle \sigma_{+},A^{(i)}\rangle$,
we get the rational  representation of them via $a_{ij}$ and $f_{ijk}$.
$\Box$

\textbf{Conclusion.}

The SchS is the Hamiltonian system on $\stackrel{\sim}{M}_{a_{\overline{ii}}}\times {\mathbb C}^{n+1}$:
$$\stackrel{\approx}{\omega}|_{\stackrel{\sim}{M}_{a_{\overline{ii}}}}-
\sum da_{ij}\wedge d\log (\lambda_i-\lambda_j)=0, \  \text{ where}
$$
$\stackrel{\sim}{M}_{a_{\overline{ii}}}$ is the submanifold $\sum A^{(i)}=0$ of
$\stackrel{\approx}{M}_{a_{\overline{ii}}}$ , of the Cartesian product of $n+1$ quadrics,
orbits $SL|\stackrel{a_{kk}}{A}\rangle'$.

The Garnier--Painlev\'e~6 system is a Hamiltonian system on
$(\stackrel{\sim}{M}_{a_{\overline{ii}}}/{SL(2,\mathbb C)} )\times
(\overline{{\mathbb C}}^{n+1}/{SL(2,\mathbb C)})$ with the same
Hamiltonians:
$$\omega - \sum da_{ij}\wedge d\log (\lambda(t_i)-\lambda(t_j))=0.$$
If $n=3$ we can set $t_1=0,t_2=1,t_3=\infty, \ t_0:=t$ and the extended phase
space is $M_{a_{\overline{ii}}}\times \overline{\mathbb{C}}\setminus\{0,1,\infty\}$.

The goal of this paper is to present the new geometrical model that
makes visible \emph{why}(how) $M_{a_{\overline{ii}}}:=\stackrel{\sim}{M}_{a_{\overline{ii}}}/{SL(2,\mathbb C)} $ is
\emph{birationally} simplectomorphic to the Cartesian product of $n-2$ quadrics $SL|\stackrel{a_{kk}}{A}\rangle '$
(to one quadric in Painlev\'e~6-case),
and not simplectomorphic. The manifold $M_{a_{\overline{ii}}}$ may be covered by
\emph{several} neighborhoods, each of them is simplectomorphic to such a product. In the case $n=3$
(Painlev\'e~6) there are three neighborhoods (quadrics).  If we add to
$M_{a_{\overline{ii}}}$ one neighborhood more, new points of which correspond
to the solutions of SchS becoming infinity in the moment $t$,
we get the so named \emph{Okamoto surface}, see \cite{Ok-s,ShT,Ba}.

In this paper we constructed the special coordinate atlas on $M_{a_{\overline{ii}}}$~,
each chart ${\mathcal{U}}( \ \dots)$ of which is isomorphic (in the Zariski topology)
to $\stackrel{\approx}{M}_{a_{\overline{ii}}\setminus\{0n-1n\}}$, to the Cartesian product of $n-2$ quadrics.
The rational simplectic map $\stackrel{\sim}{M}_{a_{\overline{ii}}}\to
\stackrel{\approx}{M}_{a_{\overline{ii}}\setminus\{0n-1n\}}$:
$$A^{(\overline{i})}
\stackrel{(\ref{eq9})\to(\ref{eq8})\to(\ref{eq5})}{
-\!\!\!-\!\!\!-\!\!\!-\!\!\!-\!\!\!
-\!\!\!-\!\!\!-\!\!\!\longrightarrow} \left(\beta_i,
q_i,
q'_i  \right)_{ i=1}^{n-2}$$
where $\beta_i=tr \ A^{(i)}\sigma_3/2$,
$q_i=tr \ A^{(i)}\sigma_{-}$, $q'_i=tr \ A^{(i)}\sigma_{+}$ and the short arrows mean the
substitutions of the corresponding formulae is presented.

The inverse map does not exist because $\tilde \pi$ is not injective,
it can be considered as the projection of the bundle
$\tilde \pi: \ \stackrel{\sim}{M}_{a_{\overline{ii}}}\to {M}_{a_{\overline{ii}}}$, its fibre is isomorphic to
$SL(2, \mathbb{C})$.
We construct \emph{the local section} of this bundle,
the rational (polynomial) map
${M}_{a_{\overline{ii}}}\supset\stackrel{\approx}{M}_{a_{\overline{ii}}\setminus\{0n-1n\}}
\longrightarrow
\stackrel{\sim}{\mathcal{U}}(\ \dots)\subset \stackrel{\sim}{M}_{a_{\overline{ii}}}$
that parameterize the fibres (formula (\ref{eq7})).

This gives the polynomial expression for the Hamiltonian's $a_{ij}=tr\ A^{(i)}A^{(j)}$ in terms
of any canonical coordinates $(p_i,q_i)_{i=1}^{n-2}$ on
$\stackrel{\approx}{M}_{a_{\overline{ii}}\setminus\{0n-1n\}}\subset{M}_{a_{\overline{ii}}}$, it is the pass
between SchS and the Garnier--Painlev\'e~6 systems announced in the Abstract.

The properties (\ref{eq2}) are evidently connected with orthonormality;
so the ${SL(2,\mathbb C)}$-invariant procedure $(\ref{eq9})\to(\ref{eq8})\to(\ref{eq5})$ is
\emph{the version of the orthonormalization of the set $A^{(\overline{i})}$ of elements of $ sl(2, \mathbb{C})$}.


\end{document}